\documentclass[a4paper]{amsart}
\usepackage{amsmath,amsthm,amsfonts,amssymb,amscd,mathrsfs}
\usepackage{psfrag,graphicx}
\usepackage{bbding}
\usepackage{epic,eepic}
\usepackage{color}
\usepackage{ebezier}

\usepackage{epsfig}

\theoremstyle{plain}
\newtheorem{main}{Theorem}

\newtheorem{Thm}{Theorem}[section]
\newtheorem{Lem}[Thm]{Lemma}

\theoremstyle{remark}
\newtheorem{Def}[Thm] {Definition}
\newtheorem{Rem}[Thm] {Remark}

\long\def\begcom#1\endcom{}

\newcommand{\length}{\operatorname{\length}}

\newcommand{\dist}{\operatorname{dist}}

\def\length{\operatorname{length}}

\def\top{\operatorname{top}}

\begin{document}

\title[Points without Physical-like Behaviour  ]
       {Topological Entropy on Points without Physical-like Behaviour}


\author[E. Catsigeras] {Eleonora Catsigeras}
\address[E. Catsigeras]{Instituto de Matem\'{a}tica y Estad\'{\i}stica \lq\lq Rafael Laguardia\rq\rq (IMERL), Facultad de Ingenier\'{\i}a,
Universidad de la Rep\'{u}blica. Uruguay
}
\email{eleonora@fing.edu.uy}

\author[X. Tian] {Xueting Tian$^1$}
\address[X. Tian]{School of Mathematical Science,  Fudan University\\Shanghai 200433, People's Republic of China}
\email{xuetingtian@fudan.edu.cn}
\urladdr{http://homepage.fudan.edu.cn/xuetingtian}

\author[E. Vargas] {Edson Vargas}
\address[E. Vargas]{Departamento de Matematica, IME-USP, Sao Paulo, Brasil}
\email{vargas@ime.usp.br}

\footnotetext[1]{X. Tian is the corresponding author.}
\maketitle

\def\abstractname{\textbf{Abstract}}

\begin{abstract}\addcontentsline{toc}{section}{\bf{English Abstract}} \small
We study a class of asymptotically entropy-expansive $C^1$ diffeomorphisms with dominated splitting on a compact manifold $M$, that satisfy the specification property. This class includes, in particular, transitive  Anosov diffeomorphisms and time-one maps of transitive Anosov flows.  We consider the nonempty set of physical-like measures that  attracts the empirical probabilities (i.e. the time averages) of Lebesgue-almost all the orbits. We define the set $I_f \cap \Gamma_f \subset M$ of irregular points without physical-like behaviour. We prove that, if not all the invariant measures of $f$ satisfy Pesin Entropy Formula (for instance in the Anosov case), then $I_f \cap \Gamma_f$ has  full topological entropy. We also obtain  this result for some class of asymptotically entropy-expansive continuous maps on a compact metric space, if the set of physical-like measures are equilibrium states with respect to some continuous potential. Finally, we prove that also the set $(M \setminus I_f) \cap \Gamma_f$ of regular points without physical-like behaviour, has full topological entropy.
\end{abstract}

\section{Introduction}

 The differentiable ergodic theory of dynamical systems is  mainly developed in the $C^{1 + \alpha}$ scenario. Relatively few results were obtained  in the $C^1$ context. In this paper we focus our attention on  $C^1$ dynamical systems with some kind of weak hyperbolicity, for which the Lebesgue measure is not necessarily invariant.

 Among the most useful concepts in the ergodic theory,  the definition of physical probability measure plays an important role.   An invariant  probability  measure $\mu$ is called physical if for a Lebesgue-positive set of initial states $x$, the time-average of any continuous function $\varphi$ along the orbit of $x$, up to time $n$, converges (when $n \rightarrow + \infty$) to the   expected value of $\varphi$ with respect to $\mu$.
Not all  systems, principally in the $C^1$ context,   possess physical measures.   This problem can be easily dodged by substituting the definition of physical measure by a weaker concept: \em physical-like \em measures, also called SRB-like or observable measures (see Definition \ref{definitionSRB-like}). Physical-like measures do always exist (see \cite{CE}).

For any $C^{1 + \alpha}$ transitive Anosov diffeomorphism $f$, the classical Pesin Theory implies the existence of a unique physical measure $\mu_f$, which is also of SRB (Sinai-Ruelle-Bowen) type. Besides, $\mu_f$ has a basin of statistical attraction with full Lebesgue measure  and satisfies Pesin Entropy Formula.  The \em typical \em points are those in the basin of statistical attraction of $\mu$.
Analogously, in   \cite{CE} was proved that for any $C^0$  system $f$,  there exists a nonempty set ${\mathcal O}_f$ composed by all the observable or physical-like measures. Besides, the basin of statistical attraction of ${\mathcal O}_f$ has full Lebesgue measure, and if $f$ is (for instance) a $C^1$ Anosov diffeomorphism, then   any measure $\mu$ in ${\mathcal O}_f$ satisfies Pesin Entropy Formula (see \cite{CCE}). In this general case, the \em typical \em points are those in the basin of statistical attraction of ${\mathcal O}_f$.

Along this paper, we will   disregard the typical orbits, and look only at  the orbits in the set of zero-Lebesgue measure that have non physical-like behaviour.  Precisely, a point $x \in M$ \em has non physical-like behaviour \em if none of the limits when $n \rightarrow \infty$ of the convergent subsequences of its time-averages, is a physical-like measure. In other words $x$ is anything but typical.

 \vspace{0cm}

We will adopt also a topologial point of view, and look at the increasing rate of the topological information quantity of  $f$; namely, its topological entropy $h_{top}(f)$.  In   \cite{Bowen},  Bowen  defined the topological entropy $h_{top}(f,E)$ restricted to an arbitrary subset $E$ of the space $M$.  Among its properties:  $h_{top}(f, E)$ increases with $E$, and $h_{top}(f, M) = h_{top}(f)$.
We say that a set $E \subset M$ \em has full topological entropy \em if $h_{top}(f, E) = h_{top}(f) $.  If so, the dynamics of $f$ restricted to $E$ produces the total increasing rate of  topological information of the system. In other words, even if one disregards the orbits whose initial states are not in $E$, the information obtained from the  sub-dynamics is, roughly speaking, the information of the whole system.

\vspace{0cm}

Our first purpose is to prove the following result that holds in the $C^1$ scenario:
\vspace{.1cm}

 \noindent{\bf Theorem 1.}  \em Let $M$ be a compact Riemannian manifold and let $f: M \mapsto M$ be a $C^1$ transitive  Anosov diffeomorphism. Then, the set of points without physical-like behaviour has full topological entropy. \em

 \vspace{.1cm}

A point $x \in M$  is \em irregular \em if the sequence of time-averages along its orbit is not convergent. The points without physical-like behaviour may be irregular or not. Besides, and principally for hyperbolic systems  that are  $C^1$ but not $C^{1+\alpha}$, the typical points  may also be irregular or not.
In \cite{Barreira-Schmeling2000} it is proved that the set of irregular points has full topological entropy. We will prove:

\vspace{.1cm}

\noindent {\bf Theorem 2.}  \em
  For any   $C^1$ transitive Anosov diffeomorphism $f$  on a compact Riemannian manifold $M$,  the set of irregular points without physical-like behaviour has full topological entropy. \em

\vspace{.1cm}
   In contrast, for any continuous map $f\colon M \mapsto M$, the set of irregular points without physical-like behaviour has zero Lebesgue measure and also zero $\mu$-measure for any $f$-invariant measure $\mu$.

   \vspace{.1cm}

We prove the above theorems using, among other tools, two well known topological properties of Anosov diffeomorphisms: the expansiveness  and the specification property. But  indeed, these conditions  are  too strong. They can be substituted by   weaker hypothesis: the asymptotical entropy-expansiveness  (see \cite{BowenEntropyExpansive,DFPV,LiaoVianaYang,PacVie}), and the $g$-almost product property for some blow-up function $g$ (see \cite{PS}). In particular expansive maps are asymptotically entropy expansive, and maps satisfying the specification property, also satisfy the $g$-almost product property. Also uniform hyperbolicity is too strong and can be weakened. We obtain:

\begin{main}
 \label{TheoremMain3} Let $f: M \mapsto M$ be a $C^1$ diffeomorphism on a compact Riemannian manifold $M$ with a dominated splitting $TM = E \oplus F$. Assume that the Lyapunov exponents  are non positive along $E$ and non negative along $F$. If $f$ is asymptotically entropy-expansive, if it satisfies the $g$-almost product property, and if not all the invariant measures satisfy Pesin Entropy Formula, then the set of irregular points without physical-like behaviour has full topological entropy.
\end{main}

Time-one diffeomorphisms of $C^1$ transitive Anosov flows  are  also particular cases for which   Theorem \ref{TheoremMain3} applies:

\vspace{.1cm}
\noindent{\bf Corollary 1.} \em
Let $f: M \mapsto M$ be the time-one map of a $C^1$ transitive Anosov flow on the compact Riemannian manifold $M$. Then, the set of irregular points without physical-like behaviour has full topological entropy. \em

\vspace{.1cm}

{\bf Proof.}  The  time-$t$ ($t\neq 0$) map of an   Anosov flow is partially hyperbolic with one-dimensional central bundle tangent to the flow. So, it has a dominated splitting $TM= E \oplus F$, where $E$ is the sum of stable and central bundles, and $F$ is the unstable bundle. Then, $f$ is far \lq\lq from  tangencies\rq\rq (see \cite{LiaoVianaYang}), and applying the known results of the theory of entropy-expansiveness (see  \cite{LiaoVianaYang} or see    \cite{DFPV,PacVie}), we deduce that $f$ is entropy-expansive. This is a stronger condition than the asymptotical entropy-expansiveness. The other hypothesis of Theorem \ref{TheoremMain3} are well known for time-$t$ maps of Anosov flows (as for  Anosov diffeomorphisms): $f$ satisfies the specification property, and the invariant measures supported on periodic orbits do not satisfy Pesin Entropy Formula.  Thus, Theorem \ref{TheoremMain3} applies. \qed

\vspace{.2 cm}

Actually, we will prove a more general result than Theorem \ref{TheoremMain3}, which we will state and prove  in Section \ref{sectionEndOfProofs} (Theorem \ref{theoremGeneral}).  Theorem \ref{TheoremMain3} is only a particular case of Theorem \ref{theoremGeneral}.

\vspace{.1cm}

We also obtain a variational principle that   holds for  continuous maps, and a generalization of Theorem \ref{TheoremMain3} in a $C^0$ setting.  Before stating this result, let us introduce some notation.
Let $f: M \mapsto M$ be a continuous map on a compact Riemannian manifold $M$. Denote by $\Delta_f  \subset M$ the set of typical points, and by $\Gamma_f \subset M$ the set of points without physical-like behaviour.  Denote by ${\mathcal P}_f$ the set of $f$-invariant probability measures, and by ${\mathcal O}_f$ the set of observable (or physical-like) measures. Denote by $I_f$ the set of irregular points. (The precise definitions of these sets are in Section \ref{sectionTopologicalDefinitions}.)

\begin{main} \label{theoremMainContinuousMaps}
 For any continuous map $f$ on a compact Riemannian manifold $M$ the following properties hold:
   \vspace{.1cm}

  \noindent{\bf   a) }
$h_{top}(f, \Delta_f) \leq \sup_{\mu \in {\mathcal O}_f} h_{\mu} (f), \ \ \ \ h_{top}(f, \Gamma_f) \leq \sup_{\mu \in {\mathcal P}_f \setminus {\mathcal O}_f} h_{\mu} (f).$

\vspace{.1cm}

 \noindent{\bf   b) }   If   $f$ satisfies the $g$-almost product property for some blow-up function $g$, then the above inequalities are equalities.

 \noindent{\bf   c) } If besides condition \em b), \em $f$ is asymptotically entropy-expansive, if   ${\mathcal O}_f$ is contained in  the set $ES_{f}(\psi )$ of equilibrium states of some continuous potencial $\psi: M \mapsto \mathbb{R}$, and if $ES_{f}(\psi ) \neq {\mathcal P}_f$, then:

  {\bf   c1) }  The set $\Gamma_f \cap I_f$ of irregular points without physical-like behaviour has full topological entropy.

{\bf   c2) }  For all $\epsilon>0$
there exists a non physical-like invariant measure $\mu$ such that
    $h_{\mu}(f) > h_{top}(f) - \epsilon$, and the set of regular points whose sequence of time-averages converge to $\mu$ has topological entropy larger than $h_{top}(f) - \epsilon$.

{\bf   c3) } The set $\Gamma_f \cap (M \setminus I_f)$ of  regular points without physical-like behaviour also has full topological entropy.

\end{main}

We remark that indeed, the diffeomorphisms satisfying the hypothesis of Theorem \ref{TheoremMain3} (in particular transitive Anosov diffeomorphisms and time-one maps of transitive Anosov flows), also satisfy the hypothesis of part c) of Theorem \ref{theoremMainContinuousMaps}. Thus,   assertions   c2) and c3) also hold for them. We will   prove in fact, a more general result than Theorem \ref{theoremMainContinuousMaps}, in Section \ref{sectionEndOfProofs} (Theorem \ref{TheoremB-2}), which contains all the results that we announced in this introduction.

To prove the theorems,  we will use  the following   main tools:  the  topological and metric properties of asymptotically entropy-expansive maps (\cite{BowenEntropyExpansive,DFPV,LiaoVianaYang,PacVie}),  the formulae of the topological entropy of saturated sets according to \cite{PS}, the basic definitions and properties of the theory of equilibrium states \cite{Keller}, and Pesin  Entropy Formula for physical-like measures of certain $C^1$ diffeomorphisms  according to \cite{CCE}.

 \vspace{.1cm}

 \noindent {\bf Organization of the paper.} Section \ref{sectionSRB-like} is a review  of definitions to make precise the statements of the theorems and their proofs.   In Section \ref{sectionLemmata} we prove two lemmas and Section \ref{sectionEndOfProofs} contains the end of the proofs.

\newpage

\section{Definitions} \label{sectionSRB-like}   \label{sectionTopologicalDefinitions}

Let $f : M \rightarrow M $ be a continuous map on a compact   manifold $M$, which does not necessarily preserve any smooth measures with respect to the Lebesgue measure.
Let $\mathcal{P}$ denote the space of all the probability measures endowed with the weak$^*$ topology, and ${\mathcal P}_f \subset {\mathcal P}$  denote the space of $f$-invariant probability measures.

\subsection{Physical-like or SRB-like measures}

\begin{Def} (Empirical probabilities or time-averages and p-omega limit.)

\label{DefinitionEmpiricalProba}
For any point $x\in M$ and for any integer number $n \geq 1$, the \em empirical probability or time-average measure \em $\Upsilon_n(x)$  of the $f$-orbit of $x$ up to time $n$, is defined by
$$\Upsilon_n(x):= \frac1n\sum_{j=0}^{n-1}\delta_{f^j(x)},$$
where $\delta_y$ is the Dirac  probability measure supported at $y\in M$. Consider the sequence  $\big\{ \Upsilon_n \big\}_{n \in \mathbb{N}^+} $ of empiric probabilities in the space ${\mathcal P}$, and define the  \em  p-omega-limit set \em  $p \omega_f(x) \subset {\mathcal P}$ as follows:
$$p \omega_f(x) := \Big\{\mu \in {\mathcal P}: \ \ \exists \ n_i \rightarrow + \infty \mbox{ such that }
 {\lim_{i \rightarrow + \infty}}\ \Upsilon_{n_i} = \mu \Big\}.$$
It is standard to check that $p \omega_f(x) \subset {\mathcal P}_f$. From \cite{DGS} we know that $p \omega_f(x)$ is always nonempty, weak$^*$-compact and connected.
 \end{Def}

\begin{Def} \label{DefinitionPhysical} (Physical or SRB measures and their basins)

We call a measure $\mu\in \mathcal P$ {\it physical or SRB} (Sinai-Ruelle-Bowen), if  the set \begin{equation}\label{BasinOfMu}A(\mu)=\{x\in M\colon \ p\omega_f(x)=\{\mu\}\}\end{equation}
 has positive Lebesgue measure. The
set $A(\mu)$ is called \em basin of statistical attraction \em of $\mu$, or in brief, basin of $\mu$ (even if $\mu$ is not physical).
\end{Def}
\begin{Rem} The above definition of physical or SRB measures is not adopted by all the authors. Some
mathematicians  require  the measure $\mu$ to be ergodic  to call it physical.
Besides, some mathematicians when studying $C^{1 + \alpha}$ systems do not define SRB as a synonym of physical measure, but take into account the property of absolute continuity on the unstable foliation. But, in the scenario of continuous systems, and even for $C^1$ systems, the unstable conditional measures can not be defined because the unstable foliation may not exist.

\end{Rem}

 \begin{Def}  {  (Physical-like measures and their $\varepsilon$-basins}, cf. \cite{CE})  
  \label{definitionSRB-like}

  Choose any metric $\dist^*$ that induces the weak$^*$ topology on the space ${\mathcal P}$ of probability measures.
A probability measure $\mu\in \mathcal{P}$ is called \em physical-like \em   (or  \em SRB-like \em or \em observable\em) if for any $\varepsilon>0$ the set \begin{equation}\label{BasinEpsilonOfMu} A_\varepsilon(\mu)=\{x\in M\colon \ \mbox{dist}^*(p\omega_f(x),\mu)<\varepsilon\}\end{equation} has positive Lebesgue measure. The set $A_\varepsilon(\mu)$ is called \em basin of $\varepsilon$-partial statistical attraction   \em of $\mu$, or in brief, $\varepsilon$-basin of $\mu$.
 We denote by $\mathcal{O}_{f}$ the set of  physical-like measures for $f$. It is standard to check that every physical-like measure is $f$-invariant and that ${\mathcal O}_f$ does not depend on the choice of the metric in ${\mathcal P}$.
\end{Def}

\begin{Def} (Basin of a compact set of probabilities)

\label{definitionA(K)} Let  $\mathcal{K}$ by a nonempty weak$^*$ compact set of probabilities.
The \em basin of statistical attraction $A(\mathcal{K})$, \em    in brief basin of ${\mathcal K}$, is defined by:

$$A(\mathcal{K}):=\{x\in M: p\omega_f(x)\subseteq \mathcal{K}\}.$$
 \end{Def}

\begin{Thm}\label{SRB-like} \label{theoremCE}  {\bf (Characterization of physical-like measures \cite{CE})}

Let $f : M \rightarrow M $ be a continuous map on a compact   manifold $M$.
 Then, the set $\mathcal{O}_{f}$ of  physical-like measures is the minimal weak$^*$ compact set
 whose basin has total Lebesgue measure.

\em In other words: $\mathcal{O}_{f}$ is nonempty, weak$^*$ compact, and  contains the limits of the convergent subsequences of the empiric probabilities for Lebesgue almost all the initial states $x \in M$. Besides, no proper subset of ${\mathcal O}_f$ has the latter three properties simultaneously.
\end{Thm}

 {\bf Proof.} See \cite{CE}.

\medskip

\subsection{Typical, irregular, and without physical-like behaviour.}\label{subsectionDelta&Gamma}

\begin{Def}
\label{definitionTypicalPoints} (The set $\Delta_f$ of typical points) \label{DefinitionTypicalPoints}

We call a point $x \in M$  \em typical \em if $pw(x)\subseteq \mathcal{O}_{f}$. Equivalently, all the convergent subsequences of empirical probabilities of $x$ converge to  physical-like measures.
We denote by $\Delta_f$ the set of typical points. Thus:
 \begin{equation}
 \label{eq:delta}\Delta_f:=\{x: pw(x)\subseteq \mathcal{O}_{f}\},.
 \end{equation}
 From Theorem \ref{SRB-like}, $\Delta_f$ has Lebesgue full measure.  \end{Def}
 \begin{Def} (The set $\Gamma_f$ of points without physical-like behaviour). \label{definitionPhysical-likeBehaviour}

 We say that a point $x \in M$ \em has not physical-like behaviour \em if $pw(x)\cap \mathcal{O}_{f}=\emptyset$. Equivalently, none of the convergent subsequences of the empirical probabilities of $x$ converge to a physical-like measure. We denote by $\Gamma_f$ the set of such points \em without physical-like behaviour. \em  Thus
  \begin{equation}
 \label{eq:without-physical}\Gamma_f=\{x: pw(x)\cap \mathcal{O}_{f}=\emptyset\}.
 \end{equation}

 Since $\Gamma_f\subseteq M\setminus \Delta_f$ we deduce that $\Gamma_f$ has   zero Lebesgue measure.

 \end{Def}
 \begin{Def} (The set $I_f$ of irregular points). \label{DefinitionIrregularPoints}

 We call a point $x \in M$ \em irregular \em if $pw(x)$ is not a singleton. Equivalently, the sequence of empirical probabilities $\Upsilon_n(x)$  does not converge.
We denote by $I_f$ the set of irregular points. Thus,
    \begin{equation}
     \label{equationI_f}
      I_f:= \big\{x \in M: \ \{\Upsilon_n(x)\}_{n \geq 1} \mbox{ does not converge}\big\}.\end{equation}
 \end{Def}

\subsection{Topological definitions. }    \label{subsectionTopologicalEntropy}

In this subsection we  list some other concepts that we will use along the proofs. Indeed, we will not formally use the mathematical conditions that impose those definitions, but  only some already known relations among them. So here, we just cite  the bibliography where the  definitions can be found.

\vspace{.1cm}

\noindent{\bf Topological entropy of a subset $E \subset M$.}  We adopt Bowen's definition of the topological entropy $h_{top}(E, f)$ of an arbitrary subset $E \subset M$, for any compact metric space $M$ and any continuous map $f$ on $M$  (see \cite{Bowen}).

 \vspace{.1cm}

\noindent {\bf Entropy-expansive  and asymptotically entropy-expansive maps.} We recall the definitions of expansive, entropy-expansive and asymptotically entropy-expansive maps in \cite{BowenEntropyExpansive,DFPV,LiaoVianaYang,PacVie}.
From those definitions, trivially every expansive homeomorphism is entropy-expansive, and every entropy expansive map (not necessarily an homeomorphism) is asymptotically entropy-expansive.

\vspace{.1cm}

\noindent {\bf Specification and $g$-almost product properties}

We adopt the definition   of the specification property of the map $f $, as for instance in \cite{DGS,Sig,Bow,Bowen2,Bowen71-trans,To2010}). We note that the original definition of specification, due to Bowen \cite{Bow}, was stronger than the specification property that we adopt here.

 We recall the definition of the blowup functions $g: \mathbb{N} \mapsto \mathbb{N}$, and of the $g$-almost product property of the map $f$, in \cite{PS}). Every continuous map $f$ that has the specification property, also has the $g$-almost product property for some blow up function $g$ (\cite[Proposition 2.1]{PS}). In other words, the $g$-almost product property is weaker than the specification property.

\subsection{Saturated sets and saturation property of the entropy}\label{subsection-saturatedsets}

Let $f: M \mapsto M$ be a continuous map on a compact metric space $M$. Recall equality (\ref{BasinOfMu}) defining the (maybe empty) basin $A(\mu)$ of  any measure $\mu$. Recall also Definition \ref{definitionA(K)} of basin $A(K)$  of any nonempty weak$^*$-compact set $K \subset {\mathcal P}_f$.

We reformulate  the definition of the saturated sets in  \cite{PS}, as follows:

\begin{Def}
\label{definitionSaturatedSet} (Saturated sets) Let $K \in {\mathcal P}_f$ be a nonempty, weak$^*$-compact and connected set of $f$-invariant probability measures. We call the (maybe empty) following set $G_K \subset M$,   \em the saturated set of $K$: \em
$$G_K=\{x\in M\colon \,pw_f(x)=K\} \subset A(K).$$
Note that $G_{\{\mu\}} = A(\mu)$ for any invariant measure $\mu$.
\end{Def}
\noindent For convenience, we introduce the following new definition inspired in the results of \cite{PS}:
 \begin{Def} (Saturation property of the entropy) \label{definitionSaturationProperty}
We say that the continuous system $f: M \mapsto M$   has the \em saturation property of the entropy, \em if  for any  nonempty, weak$^*$ compact and connected set $K \subseteq \mathcal P_f$, the following equality holds:
$$h_{top} (f,G_K)=\inf\{h_\mu (f)\colon \,\mu\in K\},$$  where $G_K \subset M$ is the saturated set of $K$,  $h_{top} (f,G_K)$ is the Bowen's topological entropy of the set $G_K$, and $h_{\mu}(f)$ is the metric entropy of $f$ with respect to the probability measure $\mu$.

Then, if $f$ has the saturation property of the entropy,   in particular for any  invariant (not necessarily ergodic) measure $\mu$, the basin of statistical attraction $A(\mu)$ satisfies \begin{equation} \label{eqn02}
h_{top}(f, A(\mu)) = h_{\mu}(f).\end{equation}
In \cite[Theorem 3]{Bowen}, Bowen proved that  equality (\ref{eqn02}) holds for any ergodic measure $\mu$, for any continuous map $f$ on a compact metric space $M$, even if $f$ does not satisfy the saturation property of the entropy. Besides, in
\cite[Theorem 1.2]{PS} it is proved  equality (\ref{eqn02}) also for any non-ergodic invariant measure, provided that   $f$ satisfies  the $g$-almost product property, even if it does not satisfy the saturation property of the entropy for other weak$^*$-compact sets $K$.

\end{Def}

\subsection{Dominated Splitting.} \label{subsectionGeneralCase} \label{SubsectionDominatedSplitting}

\begin{Def}\label{def:dominated}{ (Dominated Splitting)}

   Let $f:M\rightarrow M$ be a $C^1$ diffeomorphism on a compact Riemannian manifold $M$. Let $TM=E\oplus F$ be a  $Df$-invariant   and continuous splitting such that $dim(E)\cdot dim(F)\neq 0.$ It is called a {\it dominated splitting} if there exists $\sigma>1$ such that $$ {\|Df|_{E(x)}\|}{\|Df^{-1}_{f(x)}|_{F(f(x))}\|}\leq \sigma^{-1} \ \ \  \forall x\in M.$$
\end{Def}

\begin{Rem}
The continuity of the splitting in the latter definition is redundant (see \cite[p. 288]{BDV}). 
The classical definition of dominated splitting is  $TM=E\oplus F$ such that there exists $C>0$ and $\sigma>1$: $$ {\|Df^n|_{E(x)}\|}\cdot {\|Df^{-n}|_{F(f^n(x))}\|}\leq C \sigma^{-n}, \forall x\in M,\,\, n\geq 1.$$ It is equivalent to Definition \ref{def:dominated} (see \cite{Gour}).
\end{Rem}

\subsection{Pesin Entropy Formula} \label{subsectionPesinFormula}

\begin{Def} (Pesin Entropy Formula) \label{definitionPE_f}

Let $f: M \mapsto M$ be a $C^1$ diffeomorphism on a compact Riemannian manifold $M$ and let $\mu \in {\mathcal P}_f$. We say that $\mu$ \em satisfies Pesin Entropy Formula \em if
\begin{equation} \label{eqnPesinEntropyFormula} h_\mu(f)=  \int \sum_{\chi_i(x)\geq 0}\chi_i(x)d\mu,\end{equation}  where $h_{\mu}(f)$ is the metric entropy of $\mu$ and   $\chi_1(x)\geq\chi_2(x)\cdots \geq \chi_{\dim(M)}(x)$  denote the Lyapunov exponents of $\mu$-a.e. $x\in M.$
We denote   $$PE_f:= \big \{ \mu \in {\mathcal P}_f\colon   \ \mu \mbox{ satisfies Pesin Entropy Formula (\ref{eqnPesinEntropyFormula})}\big\}.$$
\end{Def}

\begin{Rem}
Recall that $PE_f$ is convex, because the metric entropy is affine on convex combinations of the invariant measures (\cite[Theorem 8.1]{Walter}). Besides, due to the affinity property  and   Ruelle's inequality \cite{ruelle}, either $PE_f= {\mathcal P}_f$ or the interior of $PE_f$ in ${\mathcal P}_f$ is empty. If besides $f$ is asymptotically entropy-expansive, then the entropy function $\mu \mapsto h_{\mu}(f)$ is upper semi-continuous (see \cite[Theorem 8.2]{Walter} for the expansive case, and \cite{BowenEntropyExpansive} for the entropy-expansive case; with a standard adaptation the proofs are extended to the asymptotically entropy-expansive case).
Joining the upper-continuity of the entropy function with Ruelle's inequality, it is deduced that $PE_f$ is weak$^*$-compact. We conclude that under the hypothesis of Theorem \ref{TheoremMain3}, the set $PE_f$ is (a priori maybe empty), convex,  weak$^*$-compact and with empty interior in ${\mathcal P}_f$.
\end{Rem}

\subsection{Equilibrium States} \label{subsectionEquilibriumStates} Let us return to the continuous setting and state the basic definitions and properties of the thermodynamic formalism (see for instance \cite{Keller}). Let $f$ be a continuous map on a compact manifold $M$. Fix a continuous real function $\psi:M \mapsto {\mathbb{R}}$, which is called the \em potential. \em  Consider the following real number $p_{f}(\psi)$:
$$p_{f}(\psi) :=\sup_{\mu \in {\mathcal P}_f}\Big( h_{\mu}(f) - \int \psi \, d \mu\Big).$$
The number $p_f(\psi)$ is called the \em pressure \em with respect to the potential $\psi$.

\begin{Def}
\label{definitionEquilibriumStates} The (maybe empty) set $ES_f(\psi)$ of $f$-invariant probability measures, is defined by  
$$ES_f:= \Big\{\mu \in {\mathcal P}_f\colon   h_{\mu}(f) - \int \psi \, d \mu = p_f(\psi) \Big\}.$$
The measures $\mu$ in $ES_f$ are called \em equilibrium states \em of $f$ with respect to the potential $\psi$. So    $ES_f(\psi)$ is the \em set of equilibrium states.\em
\end{Def}

\begin{Rem}
Due to the affinity property of the entropy function, $ES_f(\psi)$ is   convex and either is the whole space  ${\mathcal P}_f$, or it has empty interior in ${\mathcal P_f}$. If $f$ is asymptotically entropy-expansive, then, as said above, the entropy function is upper semi-continuous. Thus, $ES_f(\psi)$ is besides nonempty and weak$^*$-compact (see for instance \cite[Theorem 4.2.3]{Keller}).
\end{Rem}

 \section{The key lemmas.} \label{sectionLemmata}

Consider a $C^1$ diffeomorphism $f$ on a compact Riemannian manifold.

 Recall the definition of the set $\Gamma_f$ of points without physical-like behaviour, given by equality  (\ref{eq:without-physical}). Recall the notation $PE_f$ for the set of measures that satisfy Pesin Entropy Formula. Define the set $\Gamma_f^*$  by:
  \begin{equation} \label{eqnGamma_f*} \Gamma^*_f:=\{x: pw(x)\cap  {PE}_{f}=\emptyset\}.\end{equation}
  We denote by $C^0(M, \mathbb{R})$ the space of continuous real functions $\varphi: M \mapsto \mathbb{R}.$
Recall equality (\ref{equationI_f}) defining    the   set $I_f$ of irregular points. Given  $\varphi \in C^0(M, \mathbb{R})$, consider the  set $I ^{\varphi}_f$ defined by
$$I^{\varphi}_f := \Big\{x \in M: \ \frac{1}{n} \sum_{j= 0}^{n-1} \varphi\big(f^j(x)    \big) \mbox{ is not convergent }      \Big\} $$ $$ \ \ \ \ \ \ \ \ = \Big\{x \in M: \ \frac{1}{n} \sum_{j= 0}^{n-1} \int \varphi \, d \Upsilon_n(x) \mbox{ is not convergent }      \Big\}.$$
Then, due to the weak$^*$ topology, we have \begin{equation} \label{eqn01} I_f = \bigcup_{\varphi \in C^0(M,\mathbb{R})} I_f^{\varphi} .\end{equation}

\begin{Lem}\label{Main-Prop2015} \label{lemmaKEY}
Let $f : M \rightarrow M $  be a $C^1$ diffeomorphism on a compact Riemannian manifold $M$.

\vspace{.1cm}

\noindent {\bf a) } If $f$ satisfies   the saturation property of the entropy and if $\emptyset \neq PE_f\neq {\mathcal P}_f$, then:

\vspace{.1cm}

 {\bf a1) } For any $\varphi\in C^0(M, \mathbb{R})$
  such that \begin{eqnarray}\label{eq-differentintegral}
  \inf_{\omega\in\mathcal P_f}\int\varphi d\omega<\sup_{\omega\in\mathcal P_f}\int\varphi d\omega,
  \end{eqnarray}
   the set $\Gamma^*_f\cap I_f^\varphi$ carries full topological entropy.

\vspace{.1cm}

{\bf a2) }
The set $\Gamma^*_f\cap I_f$ has full topological entropy.

\vspace{.1cm}

\noindent {\bf b) } If $f$ satisfies the hypothesis of part \em a), \em and besides ${\mathcal O}_f \subset PE_f$, then the set $\Gamma_f\cap I_f$ of irregular points without physical-like behaviour  has full topological entropy.

\vspace{.1cm}

\noindent {\bf c) } If $f$ is asymptotically entropy expansive and satisfies the $g$-almost product property for some blow up function $g$, then $f$ satisfies the saturation property of the entropy.

\end{Lem}
\begin{Rem}
 \label{RemarkLemmaKEYpartC}
 We will show below that assertion c) of Lemma \ref{lemmaKEY}  holds for any continuous map $f: M \mapsto M$ on a compact metric space $M$.

\end{Rem}

\noindent {\bf Proof of Lemma \ref{lemmaKEY}.}

\vspace{.1cm}

\noindent {\bf Part a1). }  We divide the proof into two cases.

\vspace{.1cm}

{\bf Case 1.} $\sup_{\mu\in PE_f}h_\mu(f)<h_{top}(f).$\\
Fix $h$ satisfying $\sup_{\mu\in PE_f}h_\mu(f)<h<h_{top}(f).$ From the Variational Principle \cite{Walter}, there exists ${\mu_1}\in \mathcal P_f$ such that $$h_{\mu_1}(f)>h>\sup_{\mu\in PE_f}h_\mu(f).$$ Thus ${\mu_1} \in  P_f\setminus PE_f. $
  From inequality (\ref{eq-differentintegral}) we can take $\rho_0\in\mathcal P_f$ such that $$\int\varphi d\rho_0 \neq \int\varphi d{\mu_1}.$$
Take $\tau \in (0,1)$ close to 1 enough such that $\tau h_{\mu_1}(f)>h$. Construct the measure $${\mu_2}=\tau {\mu_1}+(1-\tau)\rho_0. \ \ \ \ \ \ \mbox{Then, }$$
  \begin{eqnarray}\label{eq-proof-differentintegral}
 \int\varphi d\mu_1 \neq \int\varphi d{\mu_2}, \ \ \ \ \mbox{ and }
 \end{eqnarray}
$$h_{\mu_2}(f)=\tau h_{\mu_1}(f)+(1-\tau)h_{\rho_0}(f)\geq \tau h_{\mu_1}(f)>h.$$ It follows that ${\mu_2}\in  P_f\setminus PE_f. $
\ \  Let  $$K:=\{\mu=\theta{\mu_1}+(1-\theta){\mu_2}:\,\,\theta\in[0,1]\}. \ \ \ \mbox{  Then, }$$
 $$h_\mu(f)\geq \min \{h_{\mu_1}(f),\,h_{\mu_2}(f)\}>h \ \ \mbox{ and thus }  \mu\in P_f\setminus PE_f \ \ \forall \ \mu \in K.$$ Let us prove that  $$G_K\subseteq \Gamma^*_f \cap I^\varphi_f.$$ In fact, for $x\in G_K$, we have $pw_f(x)=K$ So,   using (\ref{eq-proof-differentintegral})
  we deduce $x\in I^\varphi_f$. Besides  $pw_f(x)=K\subseteq P_f\setminus PE_f$ which implies $x\in \Gamma^*_f.$

By hypothesis, $f$ satisfies the saturation property of the entropy. Then $$h_{top}(f, \Gamma^*_f \cap I^\varphi_f)\geq h_{top}(f, G_K)=\inf_{\mu\in K}h_\mu(f)=\min \{h_{\mu_1}(f),\,h_{\mu_2}(f)\}>h.$$
Since $h$ is arbitrarily closed to $h_{\top}(f)$, the proof of Case 1 is complete.

\medskip

{\bf Case 2.} $\sup_{\mu\in PE_f}h_\mu(f)=h_{top}(f).$\\
Fix $\epsilon > 0$ and take two measures $\eta_0, \mu_1$ such that $$\eta_0\in \mathcal P_f\setminus PE_f, \ \ \  \mu_1\in PE_f \ \ \mbox{ and }  \ \  h_{\mu_1}(f)> h_{top}(f)-\epsilon.$$
 Either there exists    $$\nu_0\in \mathcal P_f\setminus PE_f \mbox{ such that } \int\varphi d\nu_0 \neq \int\varphi d{\mu_1},$$
or, $\int\varphi d\eta_0 =\int\varphi d{\mu_1}$. In this latter case,  from inequality (\ref{eq-differentintegral}) we can take $\rho_0\in\mathcal P_f$ such that $\int\varphi d\rho_0 \neq \int\varphi d{\mu_1}$.
  So, we can construct the measure ${\nu_0}=\frac12{\eta_0}+\frac12\rho_0,$ which satisfies  $$\nu_0\in\mathcal P_f\setminus PE_f  \mbox{ such that } \int\varphi d \nu_0 \neq \int\varphi d{\mu_1}.$$

 Take $t_1, t_2\in(0,1), \ t_1 \neq t_2$,  close   to 1  such that $$\min\{t_1,t_2\}h_{\mu_1}(f)> h_{top}(f)-\epsilon.$$ Construct the measures $$\mu=t_1\mu_1+(1-t_1)\nu_0, \ \ \   \nu=t_2\mu_1+(1-t_2)\nu_0. \ \ \ \mbox{ Then  }$$\begin{eqnarray}\label{eq-proof-differentintegral2222222}
 \int\varphi d\mu  \neq \int\varphi d{\nu }. \ \ \ \ \mbox{Define: } \ \ K:=\{\theta\mu+(1-\theta)\nu:\,\,\theta\in[0,1]\}.
 \end{eqnarray}
Since $\mu_1\in PE_f$ but $\nu_0 \in  {\mathcal P}_f\setminus PE_f, $ we have  $\tau\mu_1+(1-\tau)\nu_0 \in  {\mathcal P}_f\setminus PE_f$ for all $\tau \in [0,1)$. In particular, $$\mu,\nu \in  {\mathcal P}_f\setminus PE_f, \ \ \ \ \ K \subset {\mathcal P}_f \setminus PE_f.$$ Let us prove that $$G_K\subseteq \Gamma^* \cap I^\varphi_f.$$ In fact,  $pw_f(x)=K$ for all $x \in G_K$. Thus, using (\ref{eq-proof-differentintegral2222222})
  we deduce $x\in I^\varphi_f$. Besides, $pw_f(x)=K\subseteq {\mathcal P}_f\setminus PE_f $ which implies $x\in \Gamma^*_f.$
Finally, recalling that $f$ satisfies the saturation property of the entropy, we obtain $$h_{top}(f, \Gamma^*_f \cap I^\varphi_f)\geq h_{top}(f, G_K)=\inf_{\mu\in K}h_\mu(f)=\min\{h_\mu(f),h_\nu(f)\}>h_{top}(f)-\epsilon.$$
Since $\epsilon>0$ is arbitrarily small, the proof of Case 2 is complete.
We have proved part a1) of Lemma \ref{lemmaKEY}.

\vspace{.1cm}

\noindent {\bf Part a2). }
By assumption there are two different invariant measures $\mu_1\neq \mu_2.$ From Riesz Theorem, there exists a continuous function  $\varphi\in C^0(M, \mathbb{R})$ such that
    $\int\varphi d\mu_1\neq \int\varphi d\mu_2.$
    In other words, $$ \inf_{\mu\in\mathcal P_f}\int\varphi d\mu<\sup_{\mu\in\mathcal P_f}\int\varphi d\mu.$$ Applying part a1), $\Gamma^*_f\cap I_f^\varphi$ carries full topological entropy.
 Besides, $I_f^\varphi\subseteq I_f$ because of equality (\ref{eqn01}).
Therefore,   $\Gamma^*_f\cap I_f \supset \Gamma^*_f \cap I_f^{\varphi}$. Since the topological entropy of a set increases when the set increases, we deduce that $\Gamma^*_f\cap I_f$ also has full topological entropy.   We have proved part a2) of Lemma \ref{lemmaKEY}.

\vspace{.1cm}

\noindent {\bf Part b) }
Note that ${\mathcal O}_ f \subset PE_f$ implies $\Gamma_f^* \subset \Gamma_f$. Thus:
 $$ M \supset \Gamma _f\cap I_f \supset  \Gamma^*_f\cap I_f , \ \ \ \ \ \  h_{top}(f,M ) \geq h_{top}(f,\Gamma _f\cap I_f) \geq    h_{top}(f,\Gamma^* _f\cap I_f)= h_{top}(f).$$ Since $h_{top}(f,M ) = h_{top}(f)$, we deduce that the above inequalities are equalities; hence $h_{top}(f,\Gamma _f\cap I_f) = h_{top}(f)$, ending the proof of part b).

  \vspace{.1cm}

\noindent {\bf Part c). } It is an immediate corollary of the theory of saturated sets in \cite{PS}. For a seek of completeness, and also for futher use, let us explain how assertion c) follows from \cite{PS}.
From now on, and until the end of this subsection, $M$ is a compact metric space and $f: M \mapsto M$ is continuous.
Let us first recall the definitions of uniformly separated sets:
For any pair of real numbers $(\delta, \epsilon)$ such that $0 <\delta<1$ and  $\varepsilon>0$,  any pair of points $x, y \in M$ and any natural number $n \geq 1$, we say that $x,y$   are
$(\delta,n,\varepsilon)$-separated, if  $$\#\{j \in \{1,2, \ldots n\}\colon \mbox{dist}(f^jx,f^jy)>\varepsilon\} \geq \delta n.$$  A subset
$E$ is called $(\delta,n,\varepsilon)$-separated if any pair of different points of $E$ are  $(\delta,n,\varepsilon)-$separated.     Now, let us fix some notation:
 For any   $f$-invariant measure $\mu$, for any neighborhood $F $ of $\mu$ in the space of probability measures ${\mathcal P}$  (i.e. $F\subseteq \mathcal P$ is a weak$^*$-open subset  such that $\nu \in F$), and  for any natural number $n \geq 1$ we construct the following    set of points $A_n(F) \subset M$:
   $$ A_n(F):=\{x\in M|\, \Upsilon_n(x)\in F\},$$
   where  $\Upsilon_n (x)$  is the empirical probability defined in \ref{DefinitionEmpiricalProba}. We denote:
 \begin{equation} \label{eqnN(F)} N(F;\delta,n,\varepsilon) : =      \text{maximal cardinality of } \mbox{a  } (\delta, n,\varepsilon)\text{-separated subset} A_n(F). \end{equation}
From the hypothesis that $f$ is asymptotically entropy-expansive, applying \cite[Theorem 3.1]{PS}, it is deduced that there exist uniform values of $\epsilon$ and $\delta$ such that for any ergodic measure $\mu$ and any neighborhood $F \subset {\mathcal P}$ of $\mu$, if  $n$ is large enough, then  the difference $$\displaystyle h_{\mu}(f) -\frac{\log_2 N(F;\delta,n,\varepsilon)} {n}  $$ is as small as wanted. Precisely:

\vspace{.1cm}

\noindent{\bf Assertion \ref{lemmaKEY}.1. } \em For any $\eta>0,$ there exist uniform values of $\epsilon>0$ and $\delta>0$ satisfying  the following condition:  for any ergodic measure $\mu$ and any neighborhood $F\subseteq \mathcal P$ of $\mu$, there exists $n_0 \geq 1$ such that  $$\displaystyle h_{\mu}(f) -\frac{\log_2 N(F;\delta,n,\varepsilon)} {n}  < \eta \ \ \ \forall \ n \geq n_0,$$
where $N(F;\delta,n,\varepsilon)$ is defined by equality  \em (\ref{eqnN(F)}).

\vspace{.1cm}

By hypothesis $f$ has the $g$-almost product property. The Variational Principle proved in   \cite[Theorem 1.1]{PS}  states that if $f$ satisfies assertion \ref{lemmaKEY}.1 and  the $g$-almost product property, then  for any nonempty weak$^*$-compact and connected set $K \subset {\mathcal P}_f$, the saturated set $G_K \subset M$ (recall Definition \ref{definitionSaturatedSet}) satisfies $h_{\top}(f, G_K) = \inf_{\mu \in K} h_{\mu}(f)$. For convenience, we named this condition as the saturation property of the entropy (recall Definition \ref{definitionSaturationProperty}). Now, part c) of Lemma \ref{lemmaKEY} is proved.
\qed

\begin{Rem}
We could prove directly assertion a2) of Lemma \ref{lemmaKEY} without proving  assertion a1). But we preferred to prove also a1) because its statement is  stronger than  a2). In fact, $I_f $ is the   union of the noncountable familiy of sets $I_f^{\varphi}$, for all $\varphi \in C^0(M, \mathbb{R})$.
\end{Rem}

\newpage

\noindent{\bf The uniform separation property.} For convenience and further use, let us recall the following definition:

\begin{Def}  (Uniform separation property) \label{definitionUniformSeparation}
A continuous map $f$ on a compact metric space $M$ satisfies the \em uniform separation property \em if assertion \ref{lemmaKEY}.1 in the proof of Lemma \ref{lemmaKEY} holds.
\end{Def}
\begin{Rem}. \label{remarkSaturation}

(1) As proved in \cite[Theorem 3.1]{PS},  for any continuous map $f$ on a compact metric space, if $f$ is asymptotically entropy-expansive, then $f$ satisfies the uniform separation property.

(2) Besides,  in \cite[Theorem1.1]{PS}  it is proved that if a continuous map $f$ satisfies the uniform separation and the $g$-almost product properties, then $f$ has the saturation property of the entropy.
\end{Rem}

\vspace{.1cm}
 Next lemma applies to  continuous maps $f\colon M \rightarrow M $   on a compact metric space $M$. Before stating it let us fix some notation.  Consider the (necessarily nonempty and weak$^*$-compact) set ${\mathcal O}_f$. We will assume that ${\mathcal O}_f \subset ES_{f}(\psi)$ for some continuous potential $\psi$.
 Let us denote
 $$\Gamma_f(\psi) := \{ x \in M \colon pw(x) \cap ES_{f}(\psi )   \neq \emptyset\}.$$
 Recall equality (\ref{eq:without-physical}) defining the set $\Gamma_f$ of points without physical-like behaviour, and note that $$\Gamma_f(\psi) \subset \Gamma_f,$$  because ${\mathcal O}_f \subset ES_f(\psi).$
\begin{Lem} \label{lemmaKEY2}
Let $f\colon M \rightarrow M $  be a continuous map on a compact metric space $M$.
If $f$   satisfies  the saturation property of the entropy, and if besides $   {\mathcal O}_f\subset ES_f(\psi) \neq {\mathcal P}_f$ for some continuous potential $\psi: M \mapsto \mathbb{R}$, then
the set $\Gamma_f(\psi)\cap I_f$ has full topological entropy, and hence the set $\Gamma_f \cap I_f$ of irregular points without physical-like behaviour also has full topological entropy.

\end{Lem}

  {\bf Proof.} \noindent
To simplify the notation along the proof, we write $ES_f$ instead of $ES_f(\psi)$. As in the proof of Lemma \ref{lemmaKEY}, we discuss   two cases:

\vspace{.1cm}

{\bf Case 1.} $\sup_{\mu\in  {\mathcal P}_f \setminus ES_f}h_\mu(f)=h_{top}(f).$
 In this case, for any $\epsilon >0$, there exists a  measure $\mu_0 \in {\mathcal P}_f \setminus ES_f $ such that $$h_{\mu_0}(f) > h_{top}(f) - \epsilon.$$  Choose any measure $\mu_1 \in  ES_f$ and consider the segment in ${\mathcal P}_f$ composed by all the measures of the form
$$\mu_{t} := t \cdot \mu_1 + (1-t)\cdot\mu_0 , \ \ \ 0 \leq t \leq 1. $$
We have $\mu_1 \in ES_f$, and $\mu_0 \not \in ES_f$. Due to the affinity property of the entropy and of the definition of equilibrium states respect to the continuous potential $\psi$, we obtain:
$$\mu_t \not \in ES_f \ \ \ \ \ \forall \ 0 \leq t < 1. $$
Besides $\mu_t$ converges to $\mu_0$ in the strong topology when $t \rightarrow 0^+$ (i.e. $\mu_t(B) \rightarrow \mu_0(B)$ for all measurable set $B \subset M$). So, there exists $0<t_1 <1$ such that $$h_{\mu_t} > h_{top}(f) - \epsilon \ \ \ \forall \ t_1 \geq t \geq 0.$$
   Denote $$K= \{\mu_t\colon   0 \leq t \leq t_1\} \subset {\mathcal P}_f \setminus ES_f.$$
By hypothesis, $f$ has the saturation property of  the entropy. Thus, the saturated set of $K$ verifies the following equality:
$$h_{\top} (G_K) = \inf_{\mu \in K} h_{\mu}(f) = \min \big \{h_{\mu_{t_1}}(f),\  h_{\mu_{0}}(f) \big \} > h_{top}(f) - \epsilon.$$
Since $G_K = \{x \in M \colon pw(x) = K \}$, and $K$ is not a singleton, and besides it is contained in ${\mathcal P}_f \setminus ES_f$, we deduce \begin{equation} \label{eqn03} G_K \subset I_f \cap \Gamma_f (\psi),\end{equation}
Therefore $$h_{top}(I_f \cap \Gamma_f (\psi)) \geq h_{top}(G_K) > h_{top}(f) - \epsilon.$$
Since this condition holds for all $\epsilon >0$ we deduce that $$h_{top}(I_f \cap \Gamma_f (\psi)) = h_{top}(f).$$ We have proved the lemma in case 1.

\medskip

{\bf Case 2.} $\sup_{\mu\in {\mathcal P}_f \setminus ES_f}h_\mu(f)<h_{top}(f).$\\
From the Variational Principle, for all $\epsilon >0$, there exists a measure  $\mu_1$ such that
$$\mu_1 \in ES_f, \ \ \ \ h_{\mu_0}(f) > h_{top}(f) - \epsilon.$$
Take any measure $\mu_0 \not \in ES_f$, and consider the segment
in ${\mathcal P}_f$ composed by all the measures of the form
$$\mu_{t} := t \cdot \mu_1 + (1-t)\cdot\mu_0 , \ \ \ 0 \leq t \leq 1. $$
We have $\mu_1 \in ES_f$, and $\mu_0 \not \in ES_F$. Due to the affinity property of the entropy and of the definition of equilibrium states respect to the continuous potential $\psi$, we obtain:
$$\mu_t \not \in ES_f \ \ \ \ \ \forall \ 0 \leq t < 1. $$
Take $0<t_1 <1$ close to 0 such that $$h_{\mu_t} > h_{top} - \epsilon \ \ \ \forall \ 0 \leq t \leq t_1.$$
   Take a real number $t_2 \in (0, t_1). $ Denote $$K= \{\mu_t\colon   t_2 \leq t \leq t_1\} \subset {\mathcal P}_f \setminus ES_f.$$
 Thus, the saturated set of $K$ verifies:
$$h_{\top} (G_K) = \inf_{\mu \in K} h_{\mu}(f) = \min \big \{h_{\mu_{t_1}}(f),\  h_{\mu_{t_2}}(f) \big \} > h_{top}(f) - \epsilon.$$
Now, the proof finishes as in the first case (from inequality (\ref{eqn03}) to the end).
We have proved   Lemma \ref{lemmaKEY2}. \qed

\begin{Rem}
Using a similar argument to the proof of Lemma \ref{lemmaKEY}, one also can deduce that, in the  hypothesis of Lemma \ref{lemmaKEY2}, for any continuous function $\varphi: M \mapsto M$ such that
$$\inf_{\mu \in {\mathcal P}_f} \int \varphi \, d \mu < \sup_{\mu \in {\mathcal P}_f} \int \varphi \, d \mu$$
the set $I_f^{\varphi} \bigcap \Gamma_f$ has full topological entropy, which is a stronger result than the assertion that was proved in Lemma \ref{lemmaKEY2}.
\end{Rem}

\section{End of   proofs.} \label{sectionEndOfProofs}

\subsection{End of the proof of Theorem   \ref{TheoremMain3}.}

Recall Definition  \ref{definitionUniformSeparation}  of  the uniform separation property, and the statement (1) of Remark \ref{remarkSaturation}. Using them, we    reformulate Theorem \ref{TheoremMain3}  in a  more general version:

\begin{Thm}
\label{theoremGeneral}
 Let $f: M \mapsto M$ be a $C^1$ diffeomorphism on a compact Riemannian manifold $M$ such that:

 \noindent \em (1) \em $f$ has a dominated splitting $TM = E \oplus F$.

 \noindent \em (2) \em  For   every physical-like measure $\mu$, and for $\mu$ a.e. $x \in M$, the Lyapunov exponents of $x$ are non positive along $E$ and non negative along $F$.

 \noindent \em (3) \em $f$ satisfies the uniform separation property. \em

  \noindent   (4) \em   $f$ has the $g$-almost product property. \em

  \noindent   (5) \em   There exists some invariant probability measure that does not satisfy Pesin Entropy Formula.

\vspace{.1cm}

 \noindent Then, the set $I_f \cap \Gamma_f$ of irregular points without physical-like behaviour has full topological entropy.

\end{Thm}

{\bf Proof.} It is enough to check that $f$ satisfies the hypothesis of part b) of Lemma \ref{lemmaKEY}. First, $f$ has   the saturation property of the entropy because it satisfies conditions (3) and (4) (see Remark \ref{remarkSaturation}). Now, taking into account that ${\mathcal O}_f \neq \emptyset$ (see Theorem \ref{theoremCE}), it is enough to check that ${\mathcal O}_f \subset PE_f$.

To show that ${\mathcal O}_f \subset PE_f$, we first characterize the (a priori, maybe empty) set $ {\mathcal O}_ f \cap PE_f$.  For any invariant measure $\mu$, and for $\mu$-a.e. $x \in M$, denote the Lyapunov exponents of $x$ by
$\chi_1 \geq \chi_2 \geq \ldots \geq \chi_{\mbox{\footnotesize dim}(M)}.$
By hypothesis (1) and (2), for all $\mu \in {\mathcal O}_f$ and for $\mu$-a.e. $x \in M$, the following equality holds:
$$\sum_{\chi_i (x) >0} \chi_i(x) = \sum_{i= 1}^{\mbox{\footnotesize dim}(F)} \chi_i (x).$$
By Ruelle's inequality and the invariance of the sub-bundle $F$, we have
$$h_{\mu}(f) \leq \int \sum_{\chi_i   >0} \chi_i \, d \mu = \int \sum_{i= 1}^{\mbox{\footnotesize dim}(F)} \chi_i \, d \mu =\int \log |\det Df|_F| \, d \mu \ \ \ \forall \ \mu \in {\mathcal O}_f. $$
Joining the above inequality with  Definition \ref{definitionPE_f} of the set $PE_f$, we deduce that $\mu \in {\mathcal O}_f \bigcap PE_f$ if and only if
\begin{equation} \label{eqn00} h_{\mu}(f) \geq \int \log |\det Df|_F| \, d \mu. \end{equation}
 Finally, we recall \cite[Theorem 1]{CCE}, where inequality (\ref{eqn00}) is proved for any $C^1$ diffeomorphism with dominated splitting, and for any $\mu \in {\mathcal O}_f$. We conclude that ${\mathcal O}_ f \subset {\mathcal O}_f \bigcap {\mathcal P}_ f$, or equivalently  ${\mathcal O}_f \subset {\mathcal P}_ f$, as wanted. \qed

\subsection{End of the proof of Theorem \ref{theoremMainContinuousMaps}.}

{\bf a)} From the definition of   $\Delta_f$ and $\Gamma_f$ (see Subsection \ref{subsectionDelta&Gamma}), we have
 $$x \in \Delta_f \   \mbox{ if and only if } \   pw(x) \subset {\mathcal O}_f,$$
 $$x \in \Gamma_f \     \mbox{ if and only if }    pw(x) \cap {\mathcal O}_f = \emptyset \     \mbox{ if and only if }    pw (x) \subset {\mathcal P}_f \setminus {\mathcal O}_f.$$

 Define $t_0=\sup_{\mu\in \mathcal{O}_f} h_\mu(f) \geq 0.$ Thus,    $h_{\mu}(f) \leq t_0 \ \ \forall \ \mu \in pw (x), \ \ \forall \  x \in \Delta_f. $

 For any real number $t \geq 0$, define the (maybe empty) set
 $$Q(t) := \{x:\exists \mu\in pw_f(x) \,\,s.t.\,\,h_\mu(f)\leq t\}.$$
 From \cite[Theorem 2]{Bowen}:
   $h_{top}(f, Q(t))\leq t$.
So, in particular for $t= t_0$, we deduce $$\Delta_f \subset Q(t_0) \ \Rightarrow \ h_{top} (f,\Delta_f) \leq h_{top}(f,Q(t_0)) \leq t_0 = \sup_{\mu\in \mathcal{O}_f} h_\mu(f),$$
  proving the first inequality of part a) of Theorem \ref{theoremMainContinuousMaps}.
  To prove the second inequality, define
  $t_1 : = \sup_{\mu\in \mathcal P_f\setminus \mathcal{O}_f} h_\mu(f)$, and repeat the argument above,  with $ t_1$ instead of $t_0$. We have  proved part a).

\vspace{.1cm}

  \noindent{\bf b)} By hypothesis,   $f$ satisfies the $g$-almost product property.  Consider equality (\ref{BasinOfMu}) defining the basin of statistical attraction $A({\mu})$ of any (non necessarily ergodic) invariant measure $\mu$, and observe that $A(\mu)$ is the saturated set of the singleton $\{\mu\} $ (recall Definition  \ref{definitionSaturatedSet}). Now, we apply
  \cite[Theorem 1.2]{PS}, which states  that if a continuous map $f$  satisfies the $g$-almost property, then $$h_{top}\big(f,A(\mu)\big) = h_{\mu}(f)$$ for any $f$-invariant measure $\mu$.
Since for any $\mu\in \mathcal O_f$ we have $A(\mu) \subseteq \Delta_f,$ then $$h_{top}(f,\Delta_f)\geq \sup_{\mu\in \mathcal{O}_f} h_{top}(f,A(\mu))=\sup_{\mu\in \mathcal{O}_f} h_\mu(f). $$ Together with part a), we deduce that $h_{top}(f,\Delta_f)=\sup_{\mu\in \mathcal{O}_f} h_\mu(f),$ proving the first equality of part b). To prove the second equality, repeat the argument above,  with $  P_f\setminus \mathcal{O}_f$ instead of $\mathcal{O}_f$. We have  proved part b).

\vspace{.1cm}

\noindent{\bf c)}  We will prove the following more general result, substituting the hypothesis of asymptotical entropy-expansiveness by the weaker condition of the uniform separation property (recall Definition \ref{definitionUniformSeparation} and assertion (1) of Remark \ref{RemarkLemmaKEYpartC}):

\begin{Thm}
\label{TheoremB-2}
Let $f$ be a continuous map on a compact metric space $M$. If $f$ satisfies the uniform separation property and the $g$-almost product properties, and besides, if   ${\mathcal O}_f$ is contained in  the set $ES_{f}(\psi ) \neq {\mathcal P}_f $ of equilibrium states of some continuous potencial $\psi: M \mapsto \mathbb{R}$,  then  the   assertions \em c1), c2)\em and \em c3) \em of \em Theorem \ref{theoremMainContinuousMaps} \em hold.

\end{Thm}
{\bf Proof. } From the hypothesis of uniform separation and $g$-almost product, the map $f$ satisfies the saturation property of the entropy (recall Remark \ref{remarkSaturation}).

{\bf c1) } Applying Lemma \ref{lemmaKEY2}, the
  set $\Gamma_f \cap I_f$ of irregular points without physical-like behaviour has full topological entropy.

\vspace{.1cm}

 {\bf c2) }  Applying c1) and part b) of Theorem \ref{theoremMainContinuousMaps} we deduce
$$\sup _{\mu \in {\mathcal P}_f \setminus {\mathcal O}_f} h_{\mu}(f) = h_{top}({\Gamma_f}) \geq    h_{top}(\Gamma_f \cap I_f) = h_{top}(f).$$ Therefore, for all $\epsilon > 0$
there exists an invariant measure $\mu \in {\mathcal P}_f \setminus {\mathcal O}_f$ such that
    $$h_{\mu}(f) > h_{top}(f) - \epsilon.$$ Besides,   the saturation property of the entropy implies that the basin of statistical attraction $A(\mu)$ has a topological entropy equal to $h_{\mu}(f)$. Thus,
    $$h_{top}(A(\mu)) = h_{\mu}(f) > h_{top}(f)- \epsilon,$$
    as wanted.

\vspace{.1cm}

 {\bf c3) } Finally, we must prove that the set $\Gamma_f \cap ({\mathcal P} \setminus I_f)$ of  regular points without physical-like behaviour also has full topological entropy. In fact, for any $\mu \not \in {\mathcal O}_f$ any point $x$ in the basin of statistical attraction $ A(\mu) $ belongs, by definition, to the set $\Gamma_f$, and is regular because its sequence of empirical probabilities converges to $\mu$. Thus:
  $$h_{top}(A (\mu))  \leq h_{top}(\Gamma_f \cap ({M} \setminus I_f)) \ \ \ \forall \ \mu \not \in {\mathcal O}_f.$$ For all $\epsilon>0 $ we take $\mu \not \in {\mathcal O}_f $ as in part c2). We conclude:
  $$  h_{top}(\Gamma_f \cap ({ M} \setminus I_f)) > h_{top}(f) -\epsilon \ \ \ \forall \epsilon >0.$$
  Therefore, the   set $\Gamma_f \cap ({M} \setminus I_f)$ has full topological entropy, ending the proof of Theorem \ref{TheoremB-2}; hence also of   Theorem \ref{theoremMainContinuousMaps}. \qed

\section*{Acknowlegements}
 E. Catsigeras is partially supported by
Agencia Nacional de Investigaci\'{o}n e Innovaci\'{o}n (ANII) and Comisi\'{o}n Sectorial de Investigaci\'{o}n Cient\'{\i}fica (CSIC) of Univ. de la Rep\'{u}blica.   X. Tian is  supported by National Natural Science Foundation of China  (grant no. 11301088) and moreover,   he is very grateful  to FAPESP/Brazil to support  his stay
in S\~ao Paulo for a month during the preparation of this article, FAPESP grant: 2015/09133-9. E. Vargas is  partially supported by CNPq/Brazil, grant: 308171/2011-0.

\end{document}